\input amstex

\input amsppt.sty
\magnification=\magstep1
\hsize=35truecc
 \vsize=22.5truecm
\baselineskip=14truept
 \NoBlackBoxes
\def\q{\quad}
\def\qq{\qquad}
\def\mod#1{\ (\text{\rm mod}\ #1)}

\def\t{\text}
\def\qtq#1{\q\t{#1}\q}
\par\q  Preprint: September 26, 2009
\par\q
\def\f{\frac}
\def\e{\equiv}
\def\a{\alpha}
\def\b{\binom}

\def\o{\omega}

\def\({\left(}
\def\){\right)}
\def\sls#1#2{(\f{#1}{#2})}
 \def\ls#1#2{\big(\f{#1}{#2}\big)}
\def\Ls#1#2{\Big(\f{#1}{#2}\Big)}

\def\ord{\t{\rm ord}}

\let \pro=\proclaim
\let \endpro=\endproclaim

\topmatter
\title Identities and congruences for a new sequence
\endtitle
\author ZHI-Hong Sun\endauthor
\affil School of the Mathematical Sciences, Huaiyin Normal
University,
\\ Huaian, Jiangsu 223001, PR China
\\ Email: szh6174$\@$yahoo.com
\\ Homepage: http://www.hytc.edu.cn/xsjl/szh
\endaffil

 \nologo \NoRunningHeads
\abstract{Let $[x]$ be the greatest integer not exceeding $x$. In
the paper we introduce the sequence $\{U_n\}$ given by $U_0=1$ and
$U_n=-2\sum_{k=1}^{[n/2]}\b n{2k}U_{n-2k}\q(n\ge 1)$, and establish
many recursive formulas and congruences involving $\{U_n\}$.
\par\q
}
 \endabstract
 \footnote"" {The author is
supported by the Natural Sciences Foundation of China (grant No.
10971078).}
\endtopmatter
\document
\subheading{1. Introduction}

\par  The Euler numbers $\{E_n\}$ are defined by
$$E_0=1\qtq{and} E_n=-\sum_{k=1}^{[n/2]}\b n{2k}E_{n-2k}\q(n\ge 1),$$
where $[x]$ is the greatest integer not exceeding $x$. There are
many well known identities and congruences involving Euler numbers.
In the paper we introduce the sequence $\{U_n\}$ similar to Euler
numbers as below:
$$U_0=1,\q
U_n=-2\sum_{k=1}^{[n/2]}\b n{2k}U_{n-2k}\q(n\ge 1).\tag 1.1$$
Clearly $U_{2n-1}=0$ for $n\ge 1$. In Section 2
 we establish many recursive relations for $\{U_n\}$. In Section 3,
 we deduce some congruences involving $\{U_n\}$.
  As examples, for a prime $p>3$ and
$k\in\{2,4,\ldots,p-3\}$ we have
$$\sum_{x=1}^{[p/6]}\f 1{x^k}\e
\f{6^k(2^k+1)}{4(2^{k-1}+1)}\Ls p3U_{p-1-k}\mod p,$$
 where $\sls am$ is the
Legendre-Jacobi-Kronecker symbol; for a prime $p\e 1\mod 4$ we have
$$U_{\f{p-1}2}\e (1+2(-1)^{\f{p-1}4})h(-3p) \mod p,$$ where
 $h(d)$ is the class number of the form class group consisting of
 classes of primitive, integral binary quadratic forms of discriminant $d$.
 \par Let $\Bbb N$ be the set of positive integers.
   For $m\in\Bbb N$ let $\Bbb Z_m$ be the set of rational numbers
whose denominator is coprime to $m$. For a prime $p$, in [S1] the
author introduced the notion of $p$-regular functions. If $f(k)\in
\Bbb Z_p$ for $k=0,1,2,\ldots$ and $\sum_{k=0}^n\b nk(-1)^kf(k)$ $\e
0\mod{p^n}$ for all $n\in\Bbb N$, then $f$ is called a $p$-regular
function. If $f$ and $g$ are $p$-regular functions, from [S1,
Theorem 2.3] we know that $f\cdot g$ is also a $p$-regular function.
Thus all $p$-regular functions form a ring.
\par Let $p$ be an odd prime, and let $b\in\{0,2,4,\ldots\}$.
In Section 4 we show that
 $f(k)=(1-\ls
p3p^{k(p-1)+b})U_{k(p-1)+b}$ is a  $p$-regular function. Using the
properties of $p$-regular functions in [S1,S3], we deduce many
congruences for  $\{U_{2n}\}\mod {p^m}$. For example, if
$\varphi(n)$ is Euler's totient function, for $k,m\in\Bbb N$ we have
$$U_{k\varphi(p^m)+b}\e (1-\sls
p3p^{b})U_b\mod{p^m}.$$ In Section 4 we also show that $U_{2n}\e
-16n-42\mod{128}$ for $n\ge 3$.
\par  In Section 5 we show that there is a set $X$
and a map $T:\ X\rightarrow X$ such that $(-1)^nU_{2n}$
 is the number of fixed points of $T^n$.
\par In addition to the above notation, we also use throughout this
paper the following notation: $\Bbb Z\f{\q}{\q}$the set of integers,
 $\{x\}\f{\q}{\q}$the
fractional part of $x$, $\ord_pn\f{\q}{\q}$the nonnegative integer
$\a$ such that $p^{\a}\mid n$ but $p^{\a+1}\nmid n$ (that is
$p^{\a}\ \Vert\ n$), $\mu(n)\f{\q}{\q}$the M$\ddot{\t{o}}$bius
function.

 \subheading{2. Some identities involving $\{U_n\}$}
\par Let $\{U_n\}$ be defined by (1.1). Then clearly $U_n\in\Bbb Z$.
 The first few values of $U_{2n}$ are shown below:
$$\align &\ U_2=-2,\ U_4=22,\ U_6=-602,\ U_8=30742,\ U_{10}=-2523002,
\\&U_{12}=303692662,\ U_{14}=-50402079002, \ U_{16}=11030684333782.
\endalign$$

 \pro{Lemma 2.1}
We have
  $$\sum_{n=0}^{\infty}U_n\f{t^n}{n!}=\f 1{e^t+e^{-t}-1}\q\big
  (|t|<\f{\pi}3\big)$$
  and
  $$\sum_{n=0}^{\infty}(-1)^nU_{2n}\f{t^{2n}}{(2n)!}=\f 1{2\cos
  t-1}\q\big(|t|<\f{\pi}3\big).$$
 \endpro
 Proof.
By (1.1) we have
$$\align (e^t+e^{-t}-1)\Big(\sum_{n=0}^{\infty}U_n\f{t^n}{n!}\Big)
&=\Big(1+2\sum_{k=1}^{\infty}\f{t^{2k}}{(2k)!}\Big)
\Big(\sum_{m=0}^{\infty} U_m\f{t^m}{m!}\Big)
\\&=1+\sum_{n=1}^{\infty}\Big(U_n+2\sum_{k=1}^{[n/2]}\b n{2k}U_{n-2k}\Big)
\f{t^n}{n!} =1.\endalign$$ Thus,
$$\sum_{n=0}^{\infty}U_{2n}\f{t^{2n}}{(2n)!}
=\sum_{n=0}^{\infty}U_n\f{t^n}{n!}=\f 1{e^t+e^{-t}-1}.$$ Replacing
$t$ with $it$ and noting that $\t{e}^{it}+\t{e}^{-it}=2\cos t$ we
deduce the remaining result.

\par The Bernoulli numbers $\{B_n\}$ and Bernoulli polynomials
$\{B_n(x)\}$ are defined by
$$B_0=1,\ \sum_{k=0}^{n-1}\b nkB_k=0\ (n\ge 2)\qtq{and}
B_n(x)=\sum_{k=0}^n\b nkB_kx^{n-k}\ (n\ge 0).$$  The  Euler
polynomials $\{E_n(x)\}$ are
 defined by
$$\f{2\t{e}^{xt}}{\t{e}^t+1}=\sum_{n=0}^{\infty}E_n(x)\f{t^n}{n!}\
(|t|<\pi),\tag 2.1$$ which is equivalent to (see [MOS])
 $$E_n(x)+\sum_{r=0}^n\b nrE_r(x)=2x^n\
 (n\ge 0).\tag 2.2$$
  It is well known that ([MOS])
 $$\aligned E_n(x)&=\f 1{2^n}\sum_{r=0}^n\b nr(2x-1)^{n-r}E_r
\\&=\f 2{n+1}\Big(B_{n+1}(x)-2^{n+1}B_{n+1}\Big(\f
x2\Big)\Big)\\&=\f{2^{n+1}}{n+1}\Big(B_{n+1}\Ls{x+1}2-B_{n+1}\Ls
x2\Big).\endaligned\tag
 2.3$$
  In particular,
$$E_n=2^n E_n\Ls 12\qtq{and}E_n(0)=\f{2(1-2^{n+1})B_{n+1}}{n+1}.\tag 2.4$$
 It is also known that ([MOS])
 $$B_{2n+3}=0,\ B_n(1-x)=(-1)^nB_n(x)\qtq{and}E_n(1-x)=(-1)^nE_n(x).\tag 2.5$$

\pro{Lemma 2.2} For $n\in \Bbb N
$ we have
$$ E_n\Ls 13=\f {2}{n+1}\big((-2)^{n+1}-1)\big)B_{n+1}\Ls 13
= \f {2^{n+1}((-2)^{n+1}-1)}{(n+1)((-2)^n+1)}B_{n+1}\Ls 16.$$
\endpro
Proof. By (2.3) we have $E_n(\f 13)=\f 2{n+1}\big(B_{n+1}(\f
13)-2^{n+1}B_{n+1}(\f 16)\big).$ From Raabe's theorem (see
[S3,(2.9)]) we have $B_{n+1}(\f 16)+B_{n+1}(\f 16+\f
12)=2^{-n}B_{n+1}(\f 13).$ As $B_{n+1}(\f 16+\f 12)=B_{n+1}(\f
23)=(-1)^{n+1}B_{n+1}(\f 13),$ we see that
$$B_{n+1}\Ls 16=\big(2^{-n}-(-1)^{n+1}\big)B_{n+1}\Ls 13.$$
Thus,
$$\align E_n\Ls 13&=\f {2}{n+1}\Big(B_{n+1}\Ls
13-2^{n+1}B_{n+1}\Ls 16\Big)\\&=\f
{2}{n+1}\Big(1-2^{n+1}(2^{-n}-(-1)^{n+1})\Big)B_{n+1}\Ls 13\\&= \f
{2}{n+1}\cdot \f {(-2)^{n+1}-1}{2^{-n}+(-1)^n}B_{n+1}\Ls
16.\endalign$$ So the lemma is proved.

\pro{Theorem 2.1} For $n\in \Bbb N $ we have
$$U_{2n}=3^{2n}E_{2n}\Ls 13
=-2\big(2^{2n+1}+1\big)3^{2n}\f{B_{2n+1}(\f 13)}{2n+1} = -\f{
2(2^{2n+1}+1)6^{2n}}{2^{2n}+1}\cdot \f{B_{2n+1}(\f 16)}{2n+1}.$$
\endpro

Proof. Using (2.1) and Lemma 2.1 we see that
$$\align
2\sum_{n=0}^{\infty}E_{2n}\Ls 13\f{(3t)^{2n}}{(2n)!}
&=\sum_{n=0}^{\infty}E_n\Ls
13\f{(3t)^n}{n!}+\sum_{n=0}^{\infty}E_n\Ls 13\f{(-3t)^n}{n!}
\\&=\f{2e^t}{e^{3t}+1}+\f{2e^{-t}}{e^{-3t}+1}
=\f{2e^t+2e^{2t}}{e^{3t}+1}=\f{2e^t}{e^{2t}-e^t+1}
\\&=\f 2{e^t+e^{-t}-1}=2\sum_{n=0}^{\infty}U_n\f{t^n}{n!}=2\sum_{n=0}^{\infty}U_{2n}\f{t^{2n}}{(2n)!}.
\endalign$$ Thus $U_{2n}=3^{2n}E_{2n}\sls 13$. Now applying Lemma 2.2 we
deduce the remaining result.
 \pro{Theorem 2.2} For two sequences $\{a_n\}$ and $\{b_n\}$ we have
 the following inversion formula:
 $$\align &b_n=2\sum_{k=0}^{[n/2]}\b n{2k}a_{n-2k}-a_n\q(n=0,1,2,\ldots)
\\&\iff a_n=\sum_{k=0}^{[n/2]}\b
n{2k}U_{2k}b_{n-2k}\q(n=0,1,2,\ldots).\endalign$$
\endpro
Proof. It is clear that
$$\align (e^t+e^{-t}-1)\Big(\sum_{n=0}^{\infty}a_n\f{t^n}{n!}\Big)
&=\Big(-1+2\sum_{k=0}^{\infty}\f{t^{2k}}{(2k)!}\Big)
\Big(\sum_{m=0}^{\infty}a_m\f{t^m}{m!}\Big)
\\&=\sum_{n=0}^{\infty}\Big(2\sum_{k=0}^{[n/2]}\b
n{2k}a_{n-2k}-a_n\Big)\f{t^n}{n!}.\endalign$$ Thus, using Lemma 2.1
and the fact $U_{2n-1}=0$ we see that
$$\align &b_n=2\sum_{k=0}^{[n/2]}\b n{2k}a_{n-2k}-a_n\q(n=0,1,2,\ldots)
\\&\iff (e^t+e^{-t}-1)\Big(\sum_{n=0}^{\infty}a_n\f{t^n}{n!}\Big)
=\sum_{n=0}^{\infty}b_n\f{t^n}{n!}
\\&\iff \sum_{n=0}^{\infty}a_n\f{t^n}{n!}
=\Big(\sum_{n=0}^{\infty}b_n\f{t^n}{n!}\Big)
\Big(\sum_{k=0}^{\infty}U_{2k}\f{t^{2k}}{(2k)!}\Big)
\\&\iff  a_n=\sum_{k=0}^{[n/2]}\b
n{2k}U_{2k}b_{n-2k}\q(n=0,1,2,\ldots).\endalign$$ This proves the
theorem.
\pro{Theorem 2.3} Let $n$ be a nonnegative integer. For any
complex number $x$ we have
$$\align &\sum_{k=0}^{[n/2]}\b
n{2k}U_{2k}((x-1)^{n-2k}-x^{n-2k}+(x+1)^{n-2k})=x^n, \tag i
\\&\sum_{k=0}^{[n/2]}\b
n{2k}U_{2k}(x^{n-2k}+(x+3)^{n-2k})=(x+1)^n+(x+2)^n,\tag ii
\\&\sum_{k=0}^{[n/2]}\b
n{2k}U_{2k}((x+3)^{n-2k}-(x-3)^{n-2k})\tag
iii\\&\q=(x+2)^n+(x+1)^n-(x-1)^n-(x-2)^n.\endalign$$
\endpro
Proof. From the binomial theorem we see that
$$2\sum_{k=0}^{[n/2]}\b n{2k}x^{n-2k}-x^n=(x-1)^n+(x+1)^n-x^n.$$
Thus, applying Theorem 2.2 we deduce (i). Since
$$x^m-(x+1)^m+(x+2)^m+(x+1)^m-(x+2)^m+(x+3)^m=x^m+(x+3)^m,$$ from (i)
we deduce (ii). As $x^m+(x+3)^m-((x-3)^m+x^m)=(x+3)^m-(x-3)^m$, from
(ii) we deduce (iii). So the theorem is proved.

\pro{Theorem 2.4} For $n\in\Bbb N$ we have
$$\align&\sum_{k=0}^{[n/2]}\b n{2k}(2^{n-2k}-1)U_{2k}=1-U_n,\tag i
\\&\sum_{k=0}^{[(n-1)/2]}\b n{2k}6^{n-2k}U_{2k}=5^n+4^n-2^n-1,\tag ii
\\&U_{2n}=1+2^{2n}-\sum_{k=0}^n\b{2n}{2k}3^{2n-2k}U_{2k},\tag iii
\\&U_{2n}=2(-1)^n-4\sum_{k=1}^{[n/2]}\b{2n}{4k}((-4)^k-1)U_{2n-4k},\tag iv
\\&U_{2n}=4^{n-1}+\f{1+V_{2n}}4-\f
34\sum_{k=1}^{[n/3]}\b{2n}{6k}3^{6k}U_{2n-6k},\tag v\endalign$$
where $V_m$ is given by $V_0=2,\ V_1=1$ and $V_{m+1}=V_m-7V_{m-1}\
(m\ge 1).$
\endpro
Proof.  Taking $x=1$ in Theorem 2.3(i) and noting that $U_n=0$ for
odd $n$ we obtain (i). Taking $x=3$ in Theorem 2.3(iii) we deduce
(ii). Taking $x=0$ in Theorem 2.3(ii) and then replacing $n$ with
$2n$ we derive (iii). Set $i=\sqrt{-1}$. By Theorem 2.3(i) we have
$$\sum_{k=0}^n\b{2n}{2k}U_{2k}\big((i-1)^{2n-2k}-i^{2n-2k}+(i+1)^{2n-2k}\big)=i^{2n}.$$
That is,
$$\sum_{k=0}^n\b{2n}{2k}U_{2k}\big((-2i)^{n-k}-(-1)^{n-k}+(2i)^{n-k}\big)=(-1)^n.$$
Hence
$$\sum\Sb k=0\\2\mid
n-k\endSb^n\b{2n}{2k}U_{2k}\big(2^{n+1-k}(-1)^{\f{n-k}2}-1\big)+\sum\Sb
k=0\\2\nmid n-k\endSb^n\b{2n}{2k}U_{2k}=(-1)^n.$$ Therefore,
$$\sum\Sb k=0\\2\mid
n-k\endSb^n\b{2n}{2k}U_{2k}\big(2^{n+1-k}(-1)^{\f{n-k}2}-2\big)
=(-1)^n-\sum_{k=0}^n\b{2n}{2k}U_{2k}=(-1)^n-\f 12U_{2n}$$ and so
$$2\sum\Sb r=0\\2\mid
r\endSb^n\b{2n}{2r}U_{2n-2r}\big((-1)^{\f r2}2^r-1\big)=(-1)^n-\f
12U_{2n}.$$ This yields (iv).
\par Set $\omega=(-1+\sqrt{-3})/2$. From Theorem 2.3(ii) we have
$$\sum_{k=0}^n\b{2n}{2k}U_{2k}\Big
((3\o)^{2n-2k}+(3\o+3)^{2n-2k}\Big)=(3\o+1)^{2n}+(3\o+2)^{2n}.$$ It
is easily seen that $V_m=\ls{1+3\sqrt{-3}}2^m+\ls{1-
3\sqrt{-3}}2^m=(2+3\o)^m+(-1-3\o)^m$ and
$$\aligned\o^{2n-2k}+(\o+1)^{2n-2k}=\o^{2n-2k}+(\o^2)^{2n-2k}=\cases
2&\t{if}\ 3\mid n-k,\\\o+\o^2=-1&\t{if}\ 3\nmid
{n-k}.\endcases\endaligned$$ Thus $$\aligned &3\sum\Sb k=0\\3\mid
n-k\endSb^n\b{2n}{2k}3^{2n-2k}U_{2k}-\sum_{k=0}^n\b{2n}{2k}3^{2n-2k}U_{2k}
\\&=\sum_{k=0}^n\b{2n}{2k}U_{2k}\Big
((3\o)^{2n-2k}+(3\o+3)^{2n-2k}\Big)\\&=(3\o+1)^{2n}+(3\o+2)^{2n}=V_{2n}.\endaligned$$
 Hence, applying (iii) we deduce
$$\aligned &3\sum\Sb k=0\\3\mid k\endSb^ n\b{2n}{2k}3^{2k}U_{2n-2k}\\&=
3\sum\Sb k=0\\3\mid n-k\endSb^n\b{2n}{2k}3^{2n-2k}U_{2k}=
\sum_{k=0}^n\b{2n}{2k}3^{2n-2k}U_{2k}+V_{2n}\\&=1+2^{2n}-U_{2n}+V_{2n}.\endaligned$$
This yields (v). The proof is now complete. \pro{Lemma 2.3 ([MOS,
p.30])} For $n\in\Bbb N$ and $0\le x\le 1$ we have
$$E_n(x)=4\cdot\f{n!}{\pi^{n+1}}\sum_{m=0}^{\infty}\f{\sin
((2m+1)\pi x-\f{n\pi}2)}{(2m+1)^{n+1}}.$$
\endpro

\pro{Theorem 2.5} Let $n\in\Bbb N$.  Then
$$\sum_{k=0}^{\infty}\Big(\f 1{(6k+1)^{2n+1}}-\f
1{(6k+5)^{2n+1}}\Big)=(-1)^n\f{U_{2n}\cdot\pi^{2n+1}}{2\sqrt 3\cdot
3^{2n}\cdot (2n)!}.$$
\endpro
 Proof. From Lemma 2.3 and Theorem 2.1 we see that
 $$\align &(-1)^n\f{U_{2n}\cdot\pi^{2n+1}}{4\cdot 3^{2n}\cdot
{(2n)!}}\\&=(-1)^n\f{E_{2n}\ls 13\pi^{2n+1}}{4\cdot
(2n)!}=(-1)^n\sum_{m=0}^{\infty}\f{\sin
(\f{2m+1}3\pi-n\pi)}{(2m+1)^{2n+1}} \\&=\sum_{m=0}^{\infty}\f{\sin
\f{2m+1}3\pi}{(2m+1)^{2n+1}} =\f{\sqrt 3}2\sum_{k=0}^{\infty}\Big(\f
1{(6k+1)^{2n+1}}-\f 1{(6k+5)^{2n+1}}\Big).\endalign$$ This yields
the result.
 \pro{Corollary 2.1} For $n\in\Bbb N$ we have
$(-1)^nU_{2n}>0$.
\endpro
\subheading{3. Congruences involving $\{U_{2n}\}$}

\pro{Theorem 3.1} Let $p$ be a prime of the form $4k+1$. Then
$$U_{\f {p-1}2}\e\big(1+2(-1)^{\f{p-1}4}\big)h(-3p) \mod p.$$\endpro
 Proof.
From Theorem 2.1 we see that
$$\aligned U_{\f {p-1}2}&= -2(2^{\f{p+1}2}+1)
3^{\f{p-1}2}\f{B_{\f{p+1}2}(\f 13)}{\f{p+1}2}\e -4 \Big(2\Ls
2p+1\Big)\Ls 3 pB_{\f{p+1}2}\Ls 13\\&=\cases -12B_{\f{p+1}2}\ls
13\mod p&\t{if}\ p\e 1\mod{24}\\-4B_{\f{p+1}2}\ls 13\mod p&\t{if}\
p\e 5\mod{24},\\4B_{\f{p+1}2}\ls 13\mod p&\t{if}\ p\e
13\mod{24},\\12B_{\f{p+1}2}\ls 13\mod p&\t{if}\ p\e
17\mod{24}.\endcases
\endaligned$$
By [S3, Theorem 3.2(\rm i)] we have
$$h(-3p)\e\cases -4B_{\f{p+1}2}\ls 13\mod p&\t{if}\
p\e 1\mod{12},\\4B_{\f{p+1}2}\ls 13\mod p&\t{if}\ p\e
5\mod{12}\endcases$$ Now combining the above we deduce the result.
\pro{Corollary 3.1} Let $p$ be a prime of the form $4k+1$. Then
$p\nmid U_{\f{p-1}2}.$\endpro Proof. From [UW, p.40] we know that
$h(-3p)=2\sum_{a=1}^{[p/3]}\ls pa$. Thus $1\le h(-3p)<p$. Now the
result follows from Theorem 3.1.
\par For an odd prime $p$ and $a\in\Bbb Z$ with $p\nmid a$ let
$q_p(a)=(a^{p-1}-1)/p$ denote the corresponding Fermat quotient.
 \pro{Theorem 3.2}
Let $p$ be a prime greater than $5$. Then
\par $(\t{\rm i})$
$\sum\limits_{k=1}^{[p/6]}\f 1k\e -2q_p(2)-\f
32q_p(3)+p\big(q_p(2)^2+\f 34q_p(3)^2\big)-\f {5p}2\ls
p3U_{p-3}\mod{p^2},$
\par $(\t{\rm ii})$
$\sum\limits_{k=1}^{[p/3]}\f 1k\e -\f 32q_p(3)+\f 34pq_p(3)^2-
 p\ls p3U_{p-3}\mod {p^2},$
\par $(\t{\rm iii})$
$\sum\limits_{k=1}^{[2p/3]}\f{(-1)^{k-1}}k \e 9\sum\limits\Sb
k=1\\3\mid k+p\endSb^{p-1}\f 1k\e 3p\ls p3U_{p-3}\mod{p^2}.$
\par $(\t{\rm iv})$ We have
$$\align (-1)^{[\f p6]}\b{p-1}{[\f p6]}&\e 1 +p\Big(2q_p(2)+\f
32q_p(3)\Big) +p^2\Big(q_p(2)^2+3q_p(2)q_p(3) \\&\qq+\f 38q_p(3)^2
-5\Ls p3U_{p-3} \Big) \mod {p^3}\endalign$$ and
$$(-1)^{[\f p3]}\b{p-1}{[\f p3]}\e 1 +\f 32pq_p(3)+\f 38p^2q_p(3)^2
-\f {p^2}2\Ls p3U_{p-3}\mod {p^3}.$$\endpro Proof. From Theorem 2.1
and Fermat's little theorem we have
$$U_{p-3}=-\f{2(2^{p-2}+1)\cdot
6^{p-3}}{2^{p-3}+1}\cdot\f{B_{p-2}\sls 16}{p-2} \e\f 1{30}B_{p-2}\Ls
16\mod p.$$ Now applying [S4, Theorem 3.9] we deduce the result.

\pro{Theorem 3.3} Let $p>3$ be a prime and $k\in\{2,4,\ldots,p-3\}$.
Then
$$\sum_{x=1}^{[p/6]}\f 1{x^k}\e 6^k\sum\Sb x=1\\6\mid x-p\endSb^{p-1}\f 1{x^k}\e
\f{6^k(2^k+1)}{4(2^{k-1}+1)}\Ls p3U_{p-1-k}\mod p$$ and
$$\sum_{x=1}^{[p/3]}\f 1{x^k}\e 3^k\sum\Sb x=1\\3\mid x-p\endSb^{p-1}\f 1{x^k}\e
\f{6^k}{4(2^{k-1}+1)}\Ls p3U_{p-1-k}\mod p.$$
\endpro
Proof. Let $m\in\{3,6\}$. As $B_{p-k}(\f{m-1}m)=(-1)^{p-k}B_{p-k}(\f
1m)=-B_{p-k}(\f 1m)$, we see that $B_{p-k}(\{\f pm\})=\sls
p3B_{p-k}(\f 1m)$. Now putting $s=1$ and substituting $k$ with
$p-1-k$ in [S3, Corollary 2.2] we see that for
$k\in\{2,4,\ldots,p-3\}$,
$$\sum_{x=1}^{[p/m]}\f 1{x^k}\e \sum_{x=1}^{[p/m]}
x^{p-1-k}\e\f{B_{p-k}(0)-B_{p-k}(\{\f pm\})}{p-k} =-\Ls
p3\f{B_{p-k}(\f 1m)}{p-k}\mod p.$$ By [S3, (2.6)] we have
$$\sum\Sb x=1\\m\mid x-p\endSb^{p-1}\f 1{x^k}\e
\sum\Sb x=1\\m\mid x-p\endSb^{p-1} x^{p-1-k} \e  (-m)^{p-1-k}
\sum_{x=1}^{[p/m]}x^{p-1-k}\e \f 1{m^k}\sum_{x=1}^{[p/m]}\f
1{x^k}\mod p.$$
 From Theorem 2.1 we know that
$$\f{B_{p-k}(\f
16)}{p-k}=-\f{1+2^{p-1-k}}{2(2^{p-k}+1)6^{p-1-k}}U_{p-1-k}\e
-\f{1+2^{-k}}{2(2^{1-k}+1)6^{-k}}U_{p-1-k}\mod p$$ and
$$\f{B_{p-k}(\f
13)}{p-k}=-\f{U_{p-1-k}}{2\cdot 3^{p-1-k}(2^{p-k}+1)}\e
-\f{U_{p-1-k}}{2\cdot 3^{-k}(2^{1-k}+1)}\mod p.$$ Now putting all
the above together we deduce the result.

\pro{Corollary 3.2} Let $p>3$ be a prime and
$k\in\{2,4,\ldots,p-3\}$. Then
$$\sum_{x=[p/6]+1}^{[p/3]}\f 1{x^k}\e
-\f{12^k}{4(2^{k-1}+1)}\Ls p3U_{p-1-k}\mod p$$ and
$$\sum_{x=1}^{[p/3]}\f 1{x^k}\e \f 1{2^k+1}\sum_{x=1}^{[p/6]}\f 1{x^k}\e
-\f 1{2^k}\sum_{x=[p/6]+1}^{[p/3]}\f 1{x^k}\mod p.$$
\endpro\par\q
\newline{\bf Remark 3.1} For a prime $p>5$ the congruence
$\sum_{x=1}^{[p/3]}\f 1{x^2}\e \f 15\sum_{x=1}^{[p/6]}\f 1{x^2}\mod
p$ was first found by Schwindt. See [R].

\pro{Theorem 3.4} Let $p>3$ be a prime and $k\in\{2,4,\ldots,p-3\}$.
Then
$$\sum_{x=1}^{[p/3]}(-1)^{x-1}\f 1{x^k}\e -\f{3^k}2\Ls p3U_{p-1-k}\mod p$$
and
$$\sum_{x=1}^{[\f{p+3}6]}\f 1{(2x-1)^k}\e -\f{3^k}{2^{k+1}+4}\Ls p3U_{p-1-k}\mod p.$$
\endpro
 Proof.  Putting $m=3$ and $s=1$ in [S3, Corollary 2.2] and then substituting $k$
with $p-1-k$ we see that
$$\align &
E_{p-1-k}(0)-(-1)^{[\f p3]}E_{p-1-k}\Big(\Big\{\f p3\Big\}\Big)\\&\e
2(-1)^{p-1-k-1}\sum_{x=1}^{[p/3]}(-1)^xx^{p-1-k} \e
2\sum_{x=1}^{[p/3]}(-1)^{x-1}\f 1{x^k} \mod p.\endalign$$ By (2.4)
and (2.5) we have
$$ E_{p-1-k}(0)=\f{2(1-2^{p-k})B_{p-k}}{p-k}=0.$$ From (2.5) and Theorem
2.1 we have
 $$E_{p-1-k}\Big(\Big\{\f
p3\Big\}\Big)=E_{p-1-k}\Ls 13=3^{k+1-p}U_{p-1-k}\e 3^kU_{p-1-k}\mod
p.$$ Observe that $(-1)^{[\f p3]}=\sls p3$. From the above we deduce
the first part. Since
$$\sum_{x=1}^{[\f p3]}(-1)^{x-1}\f 1{x^k}=-\sum_{x=1}^{[\f p6]}\f
1{(2x)^k}+\sum_{x=1}^{[\f{p+3}6]}\f 1{(2x-1)^k},$$ applying the
first part and Theorem 3.3 we deduce the remaining result.

\pro{Corollary 3.3} Let $p$ be a prime of the form $4k+1$. Then
$$U_{\f{p-1}2}\e -2\big(2+(-1)^{\f{p-1}4}\big)
\sum_{x=1}^{[\f{p+3}6]}\Ls p{2x-1}\mod p.$$
\endpro
Proof. Taking $k=(p-1)/2$ in Theorem 3.4 and applying Euler's
criterion we obtain
$$\sum_{x=1}^{[\f{p+3}6]}\Ls {2x-1}p\e -\f{\sls 3p\sls p3}{4+2\sls
2p}U_{\f{p-1}2}=-\f 1{4+2(-1)^{\f{p-1}4}}U_{\f{p-1}2}\mod p.$$ This
yields the result.

  \subheading{4. Congruences for $U_{k(p-1)+b}\mod {p^n}$}
   \pro{Theorem 4.1} Let $n\in\Bbb N$ with $n\ge 3$, and let
$\a$ be a nonnegative integer such that $2^{\a}\mid n$. Then
$U_{2n}\e \f 23\mod{2^{\a+4}}$. Moreover,
$$U_{2n}\e \cases 48n+\f 23\mod{2^{\a+7}}&\t{if $2\mid n$,}
\\48n+22\mod{2^7}&\t{if $2\nmid n$.}\endcases$$
\endpro
Proof. From Theorem 2.4(i) we have
$$\sum_{k=0}^n\b{2n}{2k}\big(2^{2n-2k}-1\big)U_{2k}=1-U_{2n}.$$
Thus, using (1.1) we see that
$$\sum_{k=0}^n\b{2n}{2k}2^{2n-2k}U_{2k}=1+\sum_{k=0}^{n-1}\b{2n}{2k}U_{2k}=1-\f
12U_{2n}.$$ Hence
$$U_{2n}=2-2\sum_{r=0}^n\b{2n}{2r}2^{2r}U_{2n-2r}$$
and so
$$U_{2n}=\f 23\Big(1-\sum_{r=1}^n\b{2n}{2r}4^rU_{2n-2r}\Big)=
\f 23-\f{2n}3\sum_{r=1}^n\b{2n-1}{2r-1}\f{4^r}rU_{2n-2r}.\tag 4.1$$
From the definition of $U_{2n}$ we know that $2\mid U_{2m}$ for
$m\ge 1$. Thus, for $1\le r \le n$ and $n\ge 2$ we have
$\f{4^r}rU_{2n-2r}\e 0\mod 8$ and so $2n\cdot \f{4^r}rU_{2n-2r}\e
0\mod{2^{\a+4}}$. Therefore, by (4.1) we have $U_{2n}\e \f 23\mod
{2^{\a+4}}$ and hence $U_{2n}\e 6\mod{16}$ for $n\ge 2$.
\par Since $\f{4^{n-3}}n\in \Bbb Z_2$ for $n\ge 3$, we see that $\f {2n}3\cdot \f{4^n}n=
\f{2^7n}3\cdot \f{4^{n-3}}n\e 0\mod{2^{\a+7}}$. Thus, using (4.1)
and the fact $U_{2m}\e 6\mod{16}$ for $m\ge 2$ we see that for $n\ge
3$,
$$\align U_{2n}-\f 23&=-\f {2n}3\Big(4\sum_{r=1}^{n-2}
\b{2n-1}{2r-1}\f{4^{r-1}}rU_{2n-2r}-2\cdot
2^{2n-2}(2n-1)+\f{4^n}n\Big)
\\& \e-\f {2n}3\cdot4\sum_{r=1}^{n-1}
\b{2n-1}{2r-1}\f{4^{r-1}}r\cdot 6
\\&=-16n\Big(2n-1+2\b{2n-1}3+\sum_{r=3}^{n-1}\b{2n-1}{2r-1}\f{4^{r-1}}r\Big)
\\&\e -16n\Big(2n-1+2\b{2n-1}3\Big)\mod{2^{\a+7}}.\endalign$$
It is clear that
$$\align 2n-1+2\b{2n-1}3&\e 9\Big(2n-1+2\b{2n-1}3\Big)
\\&=3(2n-1)(2n+(2n-3)^2)\e 3(2n-1)(2n+1)
\\&=3(2n-1)^2+6(2n-1)
\e 4n-3\mod 8.\endalign$$ Thus,
$$U_{2n}-\f 23\e -16n(4n-3)\e 48n+32(1-(-1)^n)\mod{2^{\a+7}}.$$
This yields the result.
 \pro{Corollary 4.1} Let $n\in\Bbb N$ and $n\ge 3$. Then
$$U_{2n}\e 6\mod{16} \qtq{and}U_{2n}\e \cases 48n-42\mod{128}&\t{if $2\mid n$,}
\\-16n-42\mod{128}&\t{if $2\nmid n$.}\endcases$$
\endpro

\pro{Theorem 4.2} Let $p$ be an odd prime and
$b\in\{0,2,4,\ldots\}$. Then $f(k)=(1-\ls
p3p^{k(p-1)+b})U_{k(p-1)+b}$ is a p-regular function.
\endpro
Proof. Suppose $n\in\Bbb N$. From Theorem 2.1 and (2.3) we have
$$\align 2^{2k+b}U_{2k+b}&=2^{2k+b}\cdot 3^{2k+b}E_{2k+b}\Ls
13=3^{2k+b}\sum_{r=0}^{2k+b}\b{2k+b}r\Big(-\f 13\Big)^{2k+b-r}E_r
\\&=\sum_{r=0}^{2k+b}\b{2k+b}r(-3)^rE_r\e \sum_{r=0}^{n-1}\b{2k+b}r(-3)^rE_r
\\&=\sum_{r=0}^{n-1}(2k+b)(2k+b-1)\cdots (2k+b-r+1)\f{(-3)^r}{r!}E_r\mod{3^n}.
\endalign$$
Since $E_r\in\Bbb Z$ and $3^r/r!\in\Bbb Z_3$, there are
$a_0,a_1,\ldots,a_{n-1}\in\Bbb Z_3$ such that
$$2^{2k+b}U_{2k+b}\e a_{n-1}k^{n-1}+\cdots+a_1k+a_0\mod {3^n}\qtq{for every $k=0,1,2,\ldots$.}$$
Hence, using [S1, Theorem 2.1] we see that $2^{2k+b}U_{2k+b}$ is a
$3$-regular function. As $$\sum_{k=0}^n\b
nk(-1)^k2^{-2k-b}=2^{-b}\Big(1-\f 14\Big)^n\e 0\mod{3^n},$$ we see
that $2^{-2k-b}$ is also a $3$-regular function. Hence,   using the
above and the product theorem of $p$-regular functions (see [S1,
Theorem 2.3]) we deduce that $f(k)=U_{2k+b}$ is a $3$-regular
function. Therefore, the result is true for $p=3$.
\par Now let us consider the case $p>3$. For $x\in\Bbb Z_p$ let $\langle x\rangle_p$ be the least
nonnegative residue of $x$ modulo $p$.  Since $2\mid b$ we have
$p-1\nmid b+1$. From [S1, Theorem 3.2] we know that
$$f_1(k)=\f{B_{k(p-1)+b+1}\sls 13-p^{k(p-1)+b}B_{k(p-1)+b+1}\ls{\f
13+\langle -\f 13\rangle_p}p}{k(p-1)+b+1}$$ is a $p$-regular
function. As
$$\f{\f
13+\langle -\f 13\rangle_p}p=\cases \f{\f 13+\f{p-1}3}p=\f 13&\t{if
$p\e 1\mod 3$},\\\f{\f 13+\f{2p-1}3}p=\f 23&\t{if $p\e 2\mod
3$}\endcases$$ and $B_{k(p-1)+b+1}(\f
23)=(-1)^{k(p-1)+b+1}B_{k(p-1)+b+1}(\f 13)=-B_{k(p-1)+b+1}(\f 13)$,
we see that
$$f_1(k)=\Big(1-\Ls p3p^{k(p-1)+b}\Big)\f{B_{k(p-1)+b+1}(\f
13)}{k(p-1)+b+1}.$$ By Theorem 2.1 and the above we have
$$\align f(k)&=\Big(1-\Ls p3p^{k(p-1)+b}\Big)\cdot (-2)\big(2^{k(p-1)+b+1}+1\big)
3^{k(p-1)+b}\f{B_{k(p-1)+b+1}(\f 13)}{k(p-1)+b+1}
\\&=-2\big(2^{k(p-1)+b+1}+1\big)
3^{k(p-1)+b}f_1(k).\endalign$$ Since
$$\sum_{k=0}^n\b nk(-1)^k\big(2^{k(p-1)+b+1}+1\big)
3^{k(p-1)+b}=2\cdot 6^b(1-6^{p-1})^n+3^b(1-3^{p-1})^n\e 0\mod
{p^n},$$  using the above and the product theorem of $p$-regular
functions (see [S1, Theorem 2.3]) we deduce that $f(k)$ is a
$p$-regular function, which completes the proof.

 \par From Theorem 4.2 and [S3, Theorem 4.3 (with $t=1$ and $d=0$)] we deduce the following
result.

 \pro{Theorem 4.3} Let $p$ be an odd prime, $k,m,n\in\Bbb
N$ and $b\in\{0,2,4,\ldots\}$. Then
$$\align &\Big(1-\Ls
p3p^{k\varphi(p^m)+b}\Big)U_{k\varphi(p^m)+b} \\&\e
\sum_{r=0}^{n-1}(-1)^{n-1-r}\b{k-1-r}{n-1-r}\b
 kr\Big(1-\Ls p3p^{r\varphi(p^m)+b}\Big)U_{r\varphi(p^m)+b}\mod{p^{mn}}.
 \endalign$$
In particular, for $n=1$ we have $U_{k\varphi(p^m)+b}\e (1-\ls
p3p^b)U_b\mod{p^m}.$
\endpro

\par From Theorem 4.2 and [S1, Theorem 2.1] we deduce the following
result.

\pro{Theorem 4.4} Let $p$ be an odd prime, $n\in\Bbb N$, $p\ge n$
and $b\in\{0,2,4,\ldots\}$. Then there are unique integers
$a_0,a_1,\ldots,a_{n-1}\in\{0,\pm 1,\pm2,\ldots,$ $\pm\f{p^n-1}2\}$
such that
 $$\Big(1-\Ls p3p^{k(p-1)+b}\Big)U_{k(p-1)+b}\e
 a_{n-1}k^{n-1}+\cdots+a_1k+a_0\mod{p^n}$$ for every $k=0,1,2,\ldots.$
Moreover,  $\ord_pa_s\ge s-\ord_ps!$ for $s=0,1,\ldots,n-1$.\endpro

 \pro{Corollary 4.2} Let $k\in\Bbb N.$ Then
\par$($\rm{i})\ $U_{2k}\e -3k+1\mod{27};$
\par(\rm{ii})\ $U_{4k}\e 1250k^4+500k^3+725k^2-1205k+2\mod{3125}\ (k\ge 2);$
\par(\rm{iii})\ $U_{4k+2}\e 1250k^4-1125k^3-675k^2-52\mod {3125}.$\endpro

\par From Theorem 4.2 and [S3, Corollary 4.2(iv)] we deduce:

\pro{Theorem 4.5} Let $p$ be an odd prime, $k,m\in\Bbb N$ and
$b\in\{0,2,4,\ldots\}$. Then
$$U_{k\varphi(p^m)+b}\e (1-kp^{m-1})\Big(1-\ls
p3p^b\Big)U_b+kp^{m-1} \Big(1-\ls
p3p^{p-1+b}\Big)U_{p-1+b}\mod{p^{m+1}}.$$
\endpro

\subheading{5. $\{(-1)^nU_{2n}\}$ is realizable}
\par If $\{a_n\}_{n=1}^{\infty}$ and $\{b_n\}_{n=1}^{\infty}$ are two
sequences satisfying $a_1=b_1$ and
 $b_n+a_1b_{n-1}+\cdots+a_{n-1}b_1=na_n\ (n>1)$, following [S2] we say that
$(a_n,b_n)$ is a Newton-Euler pair. If $(a_n,b_n)$ is a Newton-Euler
pair and $a_n\in\Bbb Z$ for all $n=1,2,3,\ldots$, then we say that
$\{b_n\}$ is a Newton-Euler sequence.
 \par Let $\{b_n\}$ be a
Newton-Euler sequence. Then clearly $b_n\in\Bbb Z$ for all
$n=1,2,3,\ldots$. In [DHL], $\{-b_n\}$ is called a Newton sequence
generated by $\{-a_n\}$.
\pro{Lemma 5.1} Let
$\{b_n\}_{n=1}^{\infty}$ be a sequence of integers. Then the
following statements are equivalent:
\par $(\t{\rm i})$ $\{b_n\}$ is a
Newton-Euler sequence.
\par $(\t{\rm ii})$  $\sum_{d\mid n}\mu\ls ndb_d\e 0\mod n$ for every
$n\in\Bbb N$.
\par $(\t{\rm iii})$ For any prime $p$ and $\a,m\in\Bbb N$ with
$p\nmid m$ we have $b_{mp^{\a}}\e b_{mp^{\a-1}}\mod {p^{\a}}$.
\par $(\t{\rm iv})$ For any $n,t\in\Bbb N$ and prime $p$ with $p^t\
\Vert\ n$ we have $b_n\e b_{\f np}\mod{p^t}$.
\par $(\t{\rm v})$ There exists a sequence $\{c_n\}$ of integers such that $b_n=\sum_{d\mid n}dc_d$
for any $n\in\Bbb N$.
 \par $(\t{\rm vi})$  For any $n\in\Bbb N$ we have
     $$\sum_{k_1+2k_2+\cdots+nk_n=n}\f{b_1^{k_1}b_2^{k_2}\cdots b_n^{k_n}}
{1^{k_1}\cdot k_1!\cdot 2^{k_2}\cdot k_2!\cdots n^{k_n}\cdot
k_n!}\in\Bbb Z.$$
\par $(\t{\rm vii})$  For any $n\in\Bbb N$ we have
$$\f 1{n!}
\vmatrix b_1&b_2&b_3&\hdots&b_n
\\ -1&b_1&b_2&\hdots&b_{n-1}
\\ \q&-2&b_1&\hdots&b_{n-2}
\\&&\ddots&\ddots&\vdots
\\&&&-(n-1)&b_1\endvmatrix\in\Bbb Z.$$
\endpro
Proof.  From [A, Theorem 3] or [DHL] we know that (i), (ii) and
(iii) are equivalent. Clearly (iii) is equivalent (iv). Using the
M$\ddot{\t{o}}$bius inversion formula we see that (ii) and (v) are
equivalent. By [S2, Theorems 2.2 and 2.3], (i),(vi) and (vii) are
equivalent. So the lemma is proved.
\par\q
\par Let $\{b_n\}_{n=1}^{\infty}$ be a sequence of nonnegative
integers. If there is a set $X$ and a map $T:\ X\rightarrow X$ such
that $b_n$ is the number of fixed points of $T^n$, following [PW]
and [A] we say that $\{b_n\}$ is realizable.
\par In [PW], Puri and Ward proved that a sequence $\{b_n\}$ of
nonnegative integers is realizable if and and only if for all
$n\in\Bbb N$, $\f 1n \sum_{d\mid n}\mu(\f nd)b_d$ is a nonnegative
integer. Thus, using the M$\ddot{\t{o}}$bius inversion formula we
see that a sequence $\{b_n\}$ is realizable if and and only if there
exists a sequence $\{c_n\}$ of nonnegative integers such that
$b_n=\sum_{d\mid n}dc_d$ for any $n\in\Bbb N$.
  In [A] J. Arias de Reyna showed that $\{E_{2n}\}$ is a
Newton-Euler sequence and $\{|E_{2n}|\}$ is realizable.
\par Now we state the following result.
\pro{Theorem 5.1} $\{U_{2n}\}$ is a Newton-Euler sequence and
$\{(-1)^nU_{2n}\}$ is realizable.
\endpro
Proof. Suppose $n\in\Bbb N$ and $\a=\ord_2n$. If $2\mid n$, by
Theorem 4.1 we have $U_{2n}\e \f 23\mod {2^{\a+4}}$ and $U_n\e \f
23\mod {2^{\a+3}}$ for $n\ge 6$. Thus $U_{2n}\e \f 23\e
U_n\mod{2^{\a}}$ for $n\ge 6$. For $n=2,4$ we also have $U_{2n}\e
U_n\mod{2^{\a}}$. If $2\nmid n$, by (1.1) we have $U_{2n}\e
0=U_n\mod {2^0}$.
\par Now assume that $p$ is an odd prime divisor of $n$ and $n=p^tn_0$
 with $p\nmid n_0$.
Using Theorem 4.3 and the fact $2n_0p^{t-1}\ge t$ we see that
$$U_{2n}=U_{2n_0p^t}=U_{2n_0\varphi(p^{t})+2n_0p^{t-1}}\e
U_{2n_0p^{t-1}}\mod{p^t}.$$
\par By the above, for any prime divisor $p$ of $n$ we have
$U_{2n}\e U_{2n/p}\mod {p^t}$, where $p^t\ \Vert\ n$. Hence, it
follows from Lemma 5.1 that $\{U_{2n}\}$ is a Newton-Euler sequence.
 \par By Corollary 2.1 we have $(-1)^nU_{2n}>0$. Suppose that $p$ is
 a prime divisor of $n$ and $p^t\ \Vert\ n$.
 If $p$ is odd, then $(-1)^n=(-1)^{\f np}$. If $p=2$ and $4\mid n$,
 we have $(-1)^n=(-1)^{\f n2}$. If $p=2$ and $2\ \Vert\ n$, then
 $(-1)^n\e (-1)^{\f n2}\mod 2$. Thus, we always have $(-1)^n\e
 (-1)^{\f np}\mod {p^t}$.  By the previous argument, we also have
   $U_{2n}\e U_{2n/p}\mod {p^t}$. Therefore,
   $(-1)^nU_{2n}\e (-1)^{\f np}U_{2n/p}\mod {p^t}$.
 Hence, by Lemma 5.1 we have $\f 1n
 \sum_{d\mid n}\mu(\f nd)(-1)^dU_{2d}\in\Bbb Z$.
  Now it remains to show that
  $\sum_{d\mid n}\mu(\f nd)(-1)^dU_{2d}\ge 0$.
 \par For $m\in\Bbb N$, by Theorem 2.5 we have
 $$(-1)^mU_{2m}=\f{2\sqrt 3\cdot 3^{2m}\cdot (2m)!}{\pi^{2m+1}}
 \sum_{k=0}^{\infty}\Big(\f 1{(6k+1)^{2m+1}}- \f  1{(6k+5)^{2m+1}}\Big).$$
Since
    $$\sum_{k=0}^{\infty}\Big(\f 1{(6k+1)^{2m+1}}- \f  1{(6k+5)^{2m+1}}\Big)
    =1-\sum_{k=0}^{\infty}\Big(\f 1{(6k+5)^{2m+1}}- \f  1{(6k+7)^{2m+1}}\Big)<1$$ and
  $$\sum_{k=0}^{\infty}\Big(\f 1{(6k+1)^{2m+1}}- \f  1{(6k+5)^{2m+1}}\Big)
  >1-\f 1{5^{2m+1}}>1-\f 15=\f 45,$$
  we see that
         $$\f 45\cdot \f{2\sqrt 3\cdot 3^{2m}\cdot (2m)!}{\pi^{2m+1}} <
         (-1)^mU_{2m}<\f{2\sqrt 3\cdot 3^{2m}\cdot (2m)!}{\pi^{2m+1}}.$$
  Hence
  $$\align \sum_{d\mid n}\mu\ls nd(-1)^dU_{2d}
  &=(-1)^nU_{2n}+\sum_{d\mid n,d\le \f n2}\mu\ls nd(-1)^dU_{2d}
  \\&\ge (-1)^nU_{2n}-\sum_{1\le d\le \f n2}(-1)^dU_{2d}
  \\&>\f 45\cdot  \f{2\sqrt 3\cdot 3^{2n}\cdot (2n)!}{\pi^{2n+1}}
  -\sum_{1\le d\le \f n2} \f{2\sqrt 3\cdot 3^{2d}\cdot (2d)!}{\pi^{2d+1}}
 \\&> \f 45\cdot  \f{2\sqrt 3\cdot 3^{2n}\cdot (2n)!}{\pi^{2n+1}}
       -\sum_{d=1}^{\infty} \f{2\sqrt 3\cdot 3^{2d}\cdot n!}{\pi^{2d+1}}
 \\&=\f {8\sqrt 3}{5\pi}\cdot n!\Big\{\Ls 9{\pi^2}^n(n+1)(n+2)\cdots
 (2n)-\f 54\cdot \f{9/\pi^2}{1- 9/\pi^2}\Big\}.
 \endalign$$
 For $m\in\Bbb N$ it is clear that
 $$ \align \Ls 9{\pi^2}^{m+1}(m+2)(m+3)\cdots (2m+2)
 &=\f 9{\pi^2}(4m+2)\cdot \Ls 9{\pi^2}^m(m+1)(m+2)\cdots (2m)
 \\&>\Ls 9{\pi^2}^m (m+1)(m+2)\cdots (2m).\endalign$$
 Thus, for $n\ge 3$ we have
 $$\Ls 9{\pi^2}^n(n+1)(n+2)\cdots (2n)\ge \Ls 9{\pi^2}^3\cdot 4\cdot
 5\cdot 6>\f 54\cdot \f{9/\pi^2}{1- 9/\pi^2}$$
 and so $\sum_{d\mid n}\mu(\f nd)(-1)^dU_{2d} >0$. This inequality is
 also true for $n=1,2$. Thus,  $\{(-1)^nU_{2n}\}$ is realizable.
 This completes the proof.
\par\q
\par Let $\{a_n\}$ be defined by
$$a_1=-2\qtq{and}na_n=U_{2n}+a_1U_{2n-2}+\cdots+a_{n-1}U_2\q(n=2,3,4,\ldots).$$
By Theorem 5.1 we have $a_n\in\Bbb Z$ for all $n\in\Bbb N$. The
first few values of $a_n$ are shown below:
$$a_2=13,\ a_3=-224,\ a_4=8170,\ a_5=-522716,\ a_6=51749722,\ a_7=-7309866728.$$

\enddocument
\bye